\documentclass[12pt]{amsart}
\usepackage{amsmath}
\usepackage{amssymb}
\usepackage{amsthm}
\usepackage{hyperref}

\newtheorem{theorem}{Theorem}[section]

\theoremstyle{definition}

\newtheorem{definition}[theorem]{Definition}

\newcommand{\N}{\mathbb{N}}
\newcommand{\Z}{\mathbb{Z}}

\newcommand{\fc}{\mathcal{F}}

\DeclareMathOperator{\ord}{ord}

\DeclareMathOperator{\supp}{supp}

\DeclareMathOperator{\vo}{\mathsf{v}}

\DeclareMathOperator{\Do}{\mathsf{D}}

\DeclareMathOperator{\s}{\sigma}

\begin{document}
\title{Inverse zero-sum problems {II}}
\author{Wolfgang A. Schmid}
\address{Institut f{\"u}r Mathematik und Wissenschaftliches Rechnen,
Karl--Fran\-zens--Universit{\"a}t Graz,
Heinrichstra{\ss}e 36,
8010 Graz, Austria}
\email{wolfgang.schmid@uni-graz.at}
\thanks{Supported by the FWF (P18779-N13).}
\date{}
\maketitle

\begin{abstract}
Let $G$ be an additive finite abelian group. A sequence over $G$ is called a minimal zero-sum sequence if the sum of its terms is zero and no proper subsequence has this property. Davenport's constant of $G$ is the maximum of the lengths of the minimal zero-sum sequences over $G$. Its value is well-known for groups of rank two. We investigate the structure of minimal zero-sum sequences of maximal length for groups of rank two.
Assuming a well-supported conjecture on this problem for groups of the form $C_m \oplus C_m$, we determine the structure of these sequences
for groups of rank two. Combining our result and partial results on this conjecture, yields unconditional results for certain groups of rank two.
\end{abstract}

\section{Introduction}

In the 1960s H.~Davenport popularized the following problem, motivated by an application in algebraic number theory.
Let $G$ be an additive finite abelian group. Determine the smallest integer $\ell$ such that each sequence over $G$ of length at least $\ell$ has a non-empty subsequence the sum of whose terms equals $0 \in G$.
This integer is now called Davenport's constant of $G$, denoted $\Do(G)$.
We refer to the recent survey article \cite{gaogeroldingersurvey}, the lecture notes \cite{geroldinger_lecturenotes},  the monographs \cite{geroldingerhalterkochBOOK}, in particular Chapters 5 to 7,  and \cite{taovuBOOK}, in particular Chapter 9, for detailed information on and applications of Davenport's constant, e.g., in investigations of the arithmetic of maximal orders of algebraic number fields.

Parallel to the problem of determining Davenport's constant, a direct problem, the associated inverse problem, i.e., the problem of determining the structure of the longest sequences that do not have a subsequence with sum zero, was intensely investigated as well.
On the one hand, solutions to the inverse problem are relevant in the above mentioned applications as well, and on the other hand, inverse results for one type of group can be applied in investigations of the direct problem for other, more complicated, types of groups
(see, e.g., \cite{bhowmik3}).

In this paper, we investigate the inverse problem associated to Davenport's constant for general finite abelian groups of rank two, complementing the investigations of the first paper in this series \cite{WAS20} that focused on groups of the form $C_m^2$, i.e., the direct sum of two cyclic groups of order $m$. To put this in context, we recall that the value of Davenport's constant for groups of rank two is well-known (cf.~Theorem \ref{thm_dir} and the references there); moreover, for cyclic groups, answers to both the direct and the inverse problem are well-known (cf.~Theorems \ref{thm_dir} and \ref{thm_invcyc} and see, e.g., \cite{savchevchen07,yuan07} for refinements), whereas, for groups of rank at least three, both the direct and the inverse problem is in general wide open (see, e.g., \cite{bhowmik,WAS_c222n} for results in special cases).

For groups of the form $C_m^2$ there is a well-known and well-supported conjecture regarding the answer to the inverse problem (see Definition \ref{def_B} for details).
For groups of the form $C_2 \oplus C_{2n}$ and $C_3 \oplus C_{3n}$ the inverse problem was solved in  \cite[Section 3]{gaogeroldinger02} and \cite{chensavchev07}, respectively, and in \cite[Section 8]{gaogeroldinger99} and \cite{girard08} partial results in the general case were obtained.
Here we solve, \emph{assuming} the above mentioned conjecture for groups of the form $C_m^2$ is true, the inverse problem for general groups of rank two (see Theorem \ref{thm_new}).

In our proof, we use direct and inverse results for cyclic groups and groups of the form $C_m^2$, which we recall in Subsection \ref{sub_known}, that we combine by using the Inductive Method (cf.~\cite[Section 5.7]{geroldingerhalterkochBOOK}).

\section{Notation and terminology}
\label{sec_not}
We recall some standard notation and terminology (we follow \cite{gaogeroldingersurvey} and \cite{geroldingerhalterkochBOOK}).

We denote by $\Z$ the set of integers, and  by $\N$ and $\N_0$ the positive and non-negative integers, respectively.
For $a,b \in \Z$, we denote by $[a,b]=\{z \in \Z \colon a\le z\le b\}$, the interval of integers.
For $k \in \Z$ and $m \in \N$, we denote by $[k]_m$ the integer in $[0, m-1]$ that is congruent to $k$ modulo $m$.

Let $G$ denote an additively written finite abelian group. (Throughout, we use additive notation for abelian groups.)
For a subset $G_0 \subset G$, we denote by $\langle G_0\rangle$ the subgroup generated by $G_0$.
We call elements $e_1, \dots, e_r \in G\setminus \{0\}$ independent
if $\sum_{i=1}^rm_ie_i=0$ with $m_i \in \Z$ implies that $m_ie_i=0$ for each $i \in [1,r]$.
We call a subset of $G$ a basis if it generates $G$ and its elements are independent.
For $n \in \N$, we denote by $C_n$ a cyclic group of order $n$.
For each finite abelian group $G$, there exist uniquely determined $1< n_1 \mid \dots \mid n_r$ such that
$G\cong C_{n_1}\oplus \dots  \oplus C_{n_r}$; we refer to $r$ as the rank of $G$ and to $\exp(G)=n_r$ as the exponent of $G$.

We denote by $\fc(G)$ the, multiplicatively written,  free abelian monoid over $G$, that is,
the monoid of all formal commutative products
\[S=\prod_{g\in G} g^{\vo_g(S)}\]
with $\vo_g(S)\in \N_0$.
We call such an element $S$ a sequence over $G$.
We refer to $\vo_g(S)$ as the multiplicity of $g$ in $S$. Moreover, $\s(S)=\sum_{g \in G} \vo_g(S)g\in G$ is called the sum of $S$,
 $|S|=\sum_{g \in G} \vo_g(S)\in \N_0$ the length of $S$, and $\supp(S) = \{g \in G\colon \vo_g(S) > 0\}\subset G$ the support of $S$.

We denote the unit element of $\fc(G)$ by $1$ and call it the empty sequence.
If $T \in \fc(G)$ and $T \mid S$ (in $\fc(G)$), then we call $T$ a subsequence of $S$; we say that it is a proper subsequence if $1\neq T \neq S$.
Moreover, we denote by $T^{-1}S$ its co-divisor, i.e., the unique sequence $R$ with $RT=S$.

If $\s(S)=0$, then we call $S$ a zero-sum sequence (zss, for short), and if $\s(T)\neq 0$ for each $1 \neq T \mid S$, then we say that $S$ is zero-sum free.
We call a zss a minimal zss (mzss, for short) if it is non-empty and has no proper subsequence with sum zero.

Using the notation recalled above, the definition of Davenport's constant can be given as follows.
For a finite abelian group $G$, let $\ell\in \N$ be minimal with the property that each $S\in \fc(G)$ with $|S| \ge \ell$ has a subsequence $1\neq T\mid S$ such that $\s(T)=0$.

It is a simple and well-known fact that $\Do(G)$ is the maximal length of a mzss over $G$ and that each zero-sum free sequence of length $\Do(G)-1$ over $G$, i.e., a sequence appearing in the inverse problem associated to $\Do(G)$, is a subsequence of a mzss of length $\Do(G)$.
Since it has technical advantages, we thus in fact investigate the structure of mzss of maximal length (ml-mzss, for short) instead of zero-sum free sequences of length $\Do(G)-1$.

Each map $f: G \to G'$ between finite abelian groups extends uniquely to a monoid homomorphism $\fc(G) \to \fc(G')$, which we denote by $f$ as well.
If $f$ is a group homomorphism, then $\s(f(S))= f(\s(S))$ for each $S \in \fc(G)$.

\section{Formulation of result}
\label{sec_res}

In this section we recall the conjecture mentioned in the introduction and formulate our result.

\begin{definition}
\label{def_B}
Let $m \in \N$. The group $C_m^2$ is said to have
Property \textbf{B} if each ml-mzss equals $g^{\exp(G)-1}T$ for some $g\in C_m^2$ and $T \in \fc(C_m^2)$.
\end{definition}

Property \textbf{B} was introduced by W.~Gao and A.~Geroldinger \cite{gaogeroldinger99,gaogeroldinger03a}.
It is conjectured that for each $m\in \N$ the group $C_m^2$ has Property \textbf{B} (see the just mentioned papers and, e.g., \cite[Conjecture 4.5]{gaogeroldingersurvey}).
We recall some result on this conjecture.

By a very recent result (see \cite{gaogeroldingergryn}, and \cite{gaogeroldinger03a} for an earlier partial result) it is known that to establish Property \textbf{B} for $C_m^2$ for each $m \in \mathbb{N}$, it suffices to establish it for $C_p^2$ for each prime $p$.
Moreover, Property \textbf{B} is known to hold for $C_m^2$ for  $m \le 28$ (see \cite{bhowmik2} and \cite{gaogeroldinger03a} for $m \le 7$).
For further recent results towards establishing Property \textbf{B} see \cite{WAS15,WAS20,bhowmik2}.

As indicated in the introduction, we characterize ml-mzss for finite abelian groups of rank two, under the assumption that a certain subgroup of the group has  Property \textbf{B}.

\begin{theorem}
\label{thm_new}
Let $G$ be a finite abelian group of rank two, say, $G\cong C_m \oplus C_{mn}$ with  $m, n \in \N$ and $m \ge 2$.
The following sequences are minimal zero-sum sequences of maximal length.
\begin{enumerate}
\item \( S  = e_j^{\ord e_j-1} \prod_{i=1}^{\ord e_k} (-x_ie_j+e_k)\)
where $\{e_1,e_2\}$ is a basis of $G$ with $\ord e_2= mn$, $\{j,k\}=\{1,2\}$, and $x_i \in \N_0$ with $\sum_{i=1}^{\ord e_k}x_i \equiv -1 \pmod{\ord e_j}$.
\item
\(S=g_1^{s m -1} \prod_{i=1}^{(n+1-s)m} (-x_ig_1 + g_2)\)
where $s \in [1,n]$, $\{g_1,g_2\}$ is a generating set of $G$ with $\ord g_2=mn$ and, in case $s\neq 1$,  $mg_1=mg_2$, and
$x_i \in \N_0$ with $\sum_{i=1}^{(n+1-s)m}x_i = m -1$.
\end{enumerate}
If $C_m^2$ has Property \textbf{B}, then all minimal zero-sum sequences of maximal length over $G$ are of this form.
\end{theorem}
The case $G \cong C_m^2$, i.e.\ $n=1$, of this result is well-known and included for completeness only (see, e.g., \cite[Theorem 5.8.7]{geroldingerhalterkochBOOK}); in particular, note that (2) is redundant for $n=1$.

This result can be combined with the above mentioned results on Property \textbf{B} to yield unconditional results for special types of groups. We do not formulate these explicitly and only point out that, since $C_2^2$ and $C_3^2$ have Property \textbf{B} (cf.\ above), the  results on $C_2\oplus C_{2n}$ and $C_3 \oplus C_{3n}$ mentioned in the introduction can be obtained in this way.

\section{Proof of the result}

In this section we give the proof of Theorem \ref{thm_new}. First, we recall some results that we use in the proof. Then, we give the actual argument.

\subsection{Known results}
\label{sub_known}
The value of $\mathsf{D}(G)$ for $G$ a group of rank two, i.e., the answer to the direct problem, is well-known (see \cite{olson69_2,vanemdeboas69}).

\begin{theorem}
\label{thm_dir}
Let $m,n\in \N$. Then $\Do(C_m \oplus C_{mn})=m + mn - 1$.
\end{theorem}

Next, we recall some results on sequences over cyclic groups.
Namely, the solution to the inverse problem associated to Davenport's constant for cyclic groups, a simple special case of \cite{boveyetal},
and the Theorem of Erd{\H o}s--Ginzburg--Ziv~\cite{erdosginzetal61}

\begin{theorem}
\label{thm_invcyc}
Let $n \in \N$ and $S \in \fc(C_n)$.
\begin{enumerate}
\item $S$ is a ml-mzss if and only if $S=e^{n}$ for some $e\in C_n$ with $\langle e \rangle=C_n$.
\item If $|S|\ge 2n-1$, then there exists some $T \mid S$ with $|T|=n$ and $\s(T)=0$.
\end{enumerate}
\end{theorem}

The following result is a main tool in the proof of Theorem \ref{thm_new}. It was obtained in \cite[Proposition 4.1, Theorem 7.1]{gaogeroldinger03a};  note that, now the additional assumption in the original version (regarding the existence of zss of length $m$ and $2m$) can be dropped, since by \cite[Theorem 6.5]{gaogeroldingersurvey} it is known to be fulfilled for each $m \in \mathbb{N}$, also note that the second type of sequence requires $t \ge 3$).

\begin{theorem}
\label{thm_tm-1}
Let $m,t \in \N$ with $m\ge 2$ and $t \ge 2$. Suppose that $C_m^2$ has Property \textbf{B}. Let $S\in \fc(C_m^2)$ be a zss of length $tm-1$
that cannot be written as the product of $t$ non-empty zss.
Then for some basis $\{f_1, f_2\}$ of $C_m^2$,
\[S= f_1^{sm-1}\prod_{i=1}^{(t-s)m}(a_if_1+f_2)\]
with $s\in [1,t-1]$ and $a_i \in [0, m-1]$ where $\sum_{i=1}^{(t-s)m} a_i \equiv 1 \pmod{m}$, or
\[S= f_1^{s_1m}f_2^{s_2 m-1} (bf_1 + f_2)^{s_3 m-1} (bf_1 + 2f_2)\]
with $s_i\in \N$ such that  $s_1+s_2+s_3= t$ and $b \in [1,m-1]$ such that $\gcd\{b, m\}=1$.
\end{theorem}

\subsection{Proof of Theorem \ref{thm_new}}

We start by establishing that all the sequences are indeed ml-mzss.

Since the length of each sequence is $mn+n-1$, and by Theorem \ref{thm_dir}, it suffices to show that they are mzss.
It is readily seen that $\s(S)=0$, thus it remains to show minimality.  Let $1\neq T\mid S$ be a zss. We assert that $T=S$.
If $S$ is as given in (1), then it suffices to note that $e_j^{\ord e_j-1}$ is zero-sum free, thus $(-x_ie_j+e_k)\mid T$ for some $i\in [1, \ord e_k] $ and this implies $\prod_{i=1}^{\ord e_k} (-x_ie_j+e_k)\mid T$, which implies $S=T$.
Suppose $S$ is as given in (2). We first note that $ag_1\in \langle g_2 \rangle$ if and only if $m\mid a$. Let $v \in \N_0$ and $I \subset [1, (n+1-s)m]$ such that $T=g_1^v\prod_{i \in I}(-x_ig_1 + g_2)$. Since $\s(T)=0$ and by the above observation, it follows that $m \mid (v-\sum_{i \in I}x_i)$, say $mb=v-\sum_{i \in I}x_i$, where $b \in [0, s-1]$.  Furthermore, we get $mbg_1+|I|g_2=0$.
If $s=1$, then $b=0$, implying that $|I|=mn$ and $v=m-1$, that is $S=T$.
If $s>1$, we have $mg_1= mg_2$, thus $mn \mid |I|+mb$ and indeed $mn= |I|+mb$. Yet, $mn= |I|+mb$ implies $|I|= [1, (n+1-s)m]$ and $b=s-1$, that is $S=T$.
Thus, the sequences are mzss.

Now, we show that if $C_m^2$ has Property \textbf{B}, then each ml-mzss is of this form.

As already mentioned, the case $n=1$ is well-known (cf.~Theorem \ref{thm_tm-1}).
We thus assume $n \ge 2$, that is  $G\cong C_{m}\oplus C_{mn}$ with $m\ge 2$ and $n \ge 2$. Furthermore, let $H = \{mg \colon g \in G\} \cong C_n$ and let $\varphi : G \to G/H$ be the canonical map; we have $G/H \cong C_m^2$. We apply the Inductive Method, as in \cite[Section 8]{gaogeroldinger99}, with the exact sequence
\[0 \to H \hookrightarrow G \overset{\varphi}{\to}  G/H \to 0.\]

Let $S \in \fc(G)$ be a ml-mzss.
First, we assert that $\varphi(S)$ cannot be written as the product of $n+1$ non-empty zss, in order to apply Theorem \ref{thm_tm-1}.
Suppose this is possible, say $\varphi(S)= \prod_{i=1}^{n+1}\varphi(S_i)$ with non-empty zss $\varphi(S_i)$.
Then $\prod_{i=1}^{n+1}\s(S_i)\in \fc(H)$ has a  proper subsequence that is a zss, yielding a proper subsequence of $S$ that is a zss.

Thus, by Theorem \ref{thm_tm-1} there exists a basis $\{f_1, f_2\}$ of $C_m^2$ such that
\begin{equation}
\label{eq_struc1}
\varphi(S)= f_1^{sm-1}\prod_{i=1}^{(n+1-s)m}(a_if_1+f_2)\end{equation}
with $s\in [1,n]$, $a_i \in [0, m-1]$, and $\sum_{i=1}^{(n+1-s)m} a_i \equiv 1 \pmod{m}$ or
\begin{equation}
\label{eq_struc2}
\varphi(S)= f_1^{s_1m}f_2^{s_2 m-1} (bf_1 + f_2)^{s_3 m-1} (bf_1 + 2f_2)
\end{equation}
with $s_i\in \N$ such that  $s_1+s_2+s_3= n+1$ and $b \in [1,m-1]$ such that $\gcd(b, m)=1$.
We distinguish two cases, depending on which of the two structures $\varphi(S)$ has.

\medskip

Case 1: $\varphi(S)$ is of the form given in \eqref{eq_struc1}. Moreover, we assume the basis $\{f_1,f_2\}$ is chosen in such a way that $s$ is maximal.
Furthermore, let $\psi: G/H \to \langle f_1\rangle$ denote the projection with respect to $G/H= \langle f_1 \rangle \oplus \langle f_2 \rangle$.
Let $S=FT$ such that $\varphi(F)=  f_1^{sm-1}$ and $T= \prod_{i=1}^{(n+1-s)m} h_i$ such that  $\varphi(h_i)=a_if_1+f_2$.

We call a factorization $T= S_0S_1 \dots S_{n-s}$ admissible if  $\s(\varphi(S_i))=0$ and $|S_i|=m$ for $i \in [1, n-s]$ (then $\s(\varphi(S_0))=f_1$ and $|S_0|=m$). Since for a sequence $T'\mid T$ of length $m$ the conditions $\s(\varphi(T'))=0$ and $\s(\psi(\varphi(T')))=0$ are equivalent, the existence of admissible factorizations follows using Theorem \ref{thm_invcyc}.

Let $T= S_0S_1 \dots S_{n-s}$ be an admissible factorization such that $|\supp(S_0)|$ is maximal (among all admissible factorizations of $T$).
Moreover, let $F = F_0 F_1 \dots F_{s-1}$ with $|F_0| = m-1$ and $ |F_i| = m $ for $i \in [1,s-1]$.
Then $\s(\varphi(F_i))=0$ for $i \in [1, s-1]$, $\s(\varphi(S_i))=0$ for $i \in [1, n-s]$, and $\s(\varphi(S_0F_0))=0$.
Thus, $\s(S_0F_0)\prod_{i=1}^{s-1} \s(F_i)\prod_{i=1}^{n-s} \s(S_i)$ is a sequence over $H$, and it is a mzss.
Since its length is $n$, it follows by Theorem \ref{thm_invcyc} that there exists some generating element $e \in H$ such that this sequence is equal to $e^n$.

We show that  $|\supp(F)|=1$. We assume to the contrary that there exist distinct $g,g'\in \supp(F)$.

First, suppose $s\ge 2$.  We may assume $g\mid F_i$ and $g'\mid F_j$ for distinct $i,j\in [0, s-1]$.
Now we consider $F_i'=g^{-1}g'F_i$ and $F_j'=g'^{-1}gF_j$ and $F_k'=F_k$ for $k \notin \{i,j\}$.
As above, we get that $\s(S_0F_0')\prod_{i=1}^{s-1} \s(F_i')\prod_{i=1}^{n-s} \s(S_i)$ is a ml-mzss over $H$ and thus equal to $\bar{e}^n$ for some generating element $\bar{e}\in H$ and indeed, since at most two elements in the sequence are changed and for $n=2$ there is only one generating element of $H$, we have $e=\bar{e}$.
Thus,  $\s(F_i)=\s(F_i')= \s(F_i)+g'-g$, a contradiction.

Second, suppose $s=1$. It follows that $m \ge 3$, since for $m=2$ we have $|F|=1$.
We consider $S_0S_j$ for some $j \in [1, n-1]$. Let $S_0S_j= T'T''$ with $|T'|=|T''|=m$.
Since $\s(\varphi(T')) + \s(\varphi(T''))= f_1$ and $\s(\varphi(T')), \s(\varphi(T'')) \in \langle f_1 \rangle$ it follows that there exists some $a \in [0, m-1]$ such that $\s(\varphi(T'))= (a+1)f_1$ and  $\s(\varphi(T''))= -af_1$.
Let $F_0=F'F''$ with $|F'|= m-(a+1)$ and $|F''|=a$.
We note that $\s(T'F')\s(T''F'')\prod_{i=1, i \neq j}^{n-1} \s(S_i)$, is a ml-mzss over $H$ and again it follows that it is equal to $e^n$ (with the same element $e$ as above).
If both $F'$ and $F''$ are non-empty, we may assume $g\mid F'$ and $g'\mid F''$, to obtain a contradiction as above.
Thus, it remains to investigate whether there exists a factorization $S_0S_j= T'T''$ with $|T'|=|T''|=m$ such that
$\{\s(\varphi(T')),\s(\varphi(T''))\}\neq \{0, f_1\}$.
We observe that such a factorization exists except if $\varphi(S_0 S_j) = (bf_1+ f_2)^{2m-1}(cf_1+f_2)$ (note that $\varphi(S_0 S_j)=(bf_1+ f_2)^{2m}$ is impossible, since $\s(\varphi(S_0 S_j))\neq 0$).

Thus, if such a factorization does not exist, for each $j \in [1,n-1]$, then $\varphi(T)=(bf_1+ f_2)^{mn-1}(cf_1+f_2)$. Since $\s(\varphi(T))=f_1$, we get $cf_1=(b+1)f_1$.
Thus, with respect to the basis consisting of $\bar{f}_1=bf_1+ f_2$ and $\bar{f}_2=f_1$, we have $\varphi(S)=\bar{f}_1^{mn-1}\bar{f}_2^{m-1}(\bar{f}_1+\bar{f}_2)$, contradicting the assumption that the basis $\{f_1, f_2\}$ maximizes $s$.

Therefore, we have $|\supp(F)|=1$ and thus
\[S = g_1^{sm - 1} T\]
for some $g_1 \in G$.

First, we consider the case $s= n$.
We have $\ord g_1 = mn$ and thus $G = \langle g_1 \rangle \oplus H_2$  where $H_2 \subset G$ is a cyclic group of order $m$.
Let  $\pi:G \to H_2$ denote the projection with respect to $G=\langle g_1 \rangle \oplus H_2$.
We observe that $\pi(\prod_{i=1}^m h_i) \in \fc(H_2)$ is a mzss and consequently it is equal to $g_2^m$ for some generating element $g_2$ of $H_2$. We note that $\{g_1, g_2\}$ is a basis of $G$. Thus, $S$ is of the form given in (1).

Thus, we may assume $s<n$.
Next, we show that if $\varphi(h_j)=\varphi(h_k)$ for $j,k \in [1,(n-s+1)m]$, then $h_j=h_k$.
Since $|h_j^{-1}T|= (n-s+1 -2)m + 2m-1$ and using again the projection $\psi$ introduced above and Theorem \ref{thm_dir}, it follows that there exists
an admissible factorization  $T=S_0'S_1'\dots S_{n-s}'$ with $h_j \mid S_0'$.
Let $\ell \in [1,n-s]$ such that $h_jh_k \mid S_0'S_{\ell}'$.
Let $S_0'S_{\ell}'= T_j'T_k'$ such that $h_j \mid T_j'$, $h_k \mid T_k'$ and $|T_j'|= |T_k'|=m$.
Similarly as above, it follows that $\s(\varphi(T_j'))=(a'+1)f_1$ and  $\s(\varphi(T_k'))=-a'f_1$ for some $a'\in [0,m-1]$.
We note that $\s(T_j'g_1^{m-a'-1})\s(T_k'g_1^{a'})\s(g_1^m)^{s-1}\prod_{i=1, i \neq \ell}^{n-s} \s(S_i')$, is a ml-mzss over $H$ and thus equal to $e'^n$ for a generating element $e'\in H$.
Similarly as above, it follows that $\s(h_j^{-1}h_kT_j'g_1^{m-a'-1})=e'$. Thus, $h_j=h_k$.

Consequently, we have \[S=g_1^{sm-1}\prod_{x \in [0,m-1]} k_x^{v_x}\]
with $\varphi(k_x)=xf_1+f_2$ for $x\in [0,m-1]$ and suitable $v_x \in \N_0$.

In the following we show that $S$ is of the form given in (2) or $\ord g_1=m$.
At the end we show that if $\ord g_1=m$, then $S$ is of the form given in (1).
We start with the following assertion.

\medskip
\noindent
\textbf{Assertion:}
Let $T= \bar{S_0}\bar{S_1} \dots \bar{S}_{n-s}$ be an admissible factorization.
Let $k_x \mid \bar{S}_0$ and let $k_y\mid \bar{S}_i$ for some $i \in[1, n-s]$. If $x< y$, then $k_y-k_x= (y-x)g_1$ and if
$x> y$, then $k_y - k_x=  (y-x)g_1 +mg_1$.

\noindent
\emph{Proof of Assertion:}
We note that $\s(\varphi(\bar{S}_0k_x^{-1}k_y))=(-x+y+1)f_1$ and  $\s(\varphi(\bar{S}_ik_y^{-1}k_x))=(-y+x)f_1$. Thus,
we have $\s(\varphi(\bar{S}_0k_x^{-1}k_yg_1^{[x-y-1]_m}))=0= \s(\varphi(\bar{S}_ik_y^{-1}k_xg_1^{[y-x]_m}))$.
We observe that $[x-y-1]_m+[y-x]_m= m-1$. Thus, similarly as above, $\s(\bar{S}_i)= \s(\bar{S}_ik_y^{-1}k_xg_1^{[y-x]_m})$. Thus, $k_y= k_x+[y-x]_m g_1$.
Consequently, if $x< y$, then $k_y-k_x= (y-x)g_1$ and if
$x> y$, then $k_y -k_x=  (y-x)g_1 +mg_1$, proving the assertion.

\medskip

First, we show that $\supp(S_0^{-1}T)\subset \supp(S_0)$ or $\ord g_1= m$.
We assume that there exists some $i \in [1, n-s]$ and some $k_t \mid S_i$  such that $k_t \nmid S_0$ and show that this implies $\ord g_1=m$.

The sequence $k_t^{-1}S_iS_0$ has length $2m-1$.
Thus, as above, there exists a subsequence $S_i''\mid k_t^{-1}S_iS_0$ such that $\s(\varphi(S_i''))=0$ and  $|S_i''|=m$. Let $S_0''= S_i''^{-1}S_iS_0$ and $S_j''=S_j$ for $j \notin \{0, i\}$.
We get that $S_0''S_1''\dots S_{n-s}''$ is an admissible factorization of $T$.
Since $k_t\mid S_0''$ and $k_t \nmid S_0$ and $|\supp(S_0)|$ is maximal (by assumption), there exists some $k_u\mid S_0$ such that $k_u \nmid S_0''$ and thus $k_u \mid S_i''$. Clearly $k_u \neq k_t$ and thus $t\neq u$.

We apply the Assertion twice.
First, to $k_u \mid S_0$ and $k_t \mid S_i$.
If $u< t$, then $k_t - k_u = (t-u)g_1$ and if
$u> t$, then $k_t - k_u = (t-u)g_1 + mg_1$.
Second, to $k_t \mid S_0''$ and $k_u \mid S_i''$.
If $t< u$, then $k_u - k_t = (u-t)g_1$ and if
$t> u$, then $k_u - k_t = (u-t)g_1 + mg_1$.

Thus, if $u< t$, then $k_t-k_u= (t-u)g_1$ and $k_u - k_t=  (u-t)g_1 +mg_1$. Adding these two equations, we get $mg_1=0$ and $\ord g_1= m$.
And, if $u> t$, then $k_t - k_u=  (t-u)g_1 +mg_1$  and  $k_u-k_t= (u-t)g_1$, again yielding  $\ord g_1= m$.

Second, we show that $|\supp(S_0^{-1}T)|=1$ or $\ord g_1=m$.
We assume that $|\supp(S_0^{-1}T)|\ge 2$, say it contains elements $k_u, k_t$ with $t> u$.
Let $i, j \in [1, n-s]$, not necessarily distinct, such that $k_u\mid S_i$ and $k_t \mid S_j$.
By the above argument we may assume that $\supp(S_0^{-1}T) \subset \supp(S_0)$.
We apply the Assertion with $k_t \mid S_0$ and $k_u\mid S_i$, to obtain $k_u -k_t=  (u-t)g_1 +mg_1$.
And, we apply the Assertion with $k_u \mid S_0$ and $k_t \mid S_j$ to obtain $k_t-k_u= (t-u)g_1$. Thus, we obtain $mg_1=0$.

Consequently, we have $T=k_u^{(n-s)m}S_0$ and $k_u \mid S_0$ for some $u\in [0, m-1]$ or $\ord g_1 = m$.

We assume that $T=k_u^{(n-s)m}S_0$ and $k_u \mid S_0$ for some $u\in [0, m-1]$.
Since $n-s \ge 1$ and $\s(k_u^m)=e$, it follows that $\ord k_u= mn$. Let $f_2'=\varphi(k_u)= uf_1+f_2$. It follows that $\{f_1, f_2'\}$ is a basis of $G/H$.

If, for $x\in [0, m-1]$, an element $h \in \supp (T)=\supp(S_0)$ exists with $\varphi(h)= -xf_1+f_2'$ (as shown above there is at most one such element), then we denote it by $k_{-x}'$. In particular, $k_u=k_0'$.

For each $k_{-x}'\in \supp(S_0)$, similarly as above, $\s(k_0'^{m}) = \s(k_0'^{m-1}k_{-x}'g_1^x)$. Thus, we have $k_0'= k_{-x}' + x g_1$.
Let $x_i \in [0,m-1]$ such that $S_0=\prod_{i=1}^m k_{-x_i}'$. We know that $\sum_{i=1}^m (-x_if_1)=f_1$, i.e., $\sum_{i=1}^m x_i\equiv m-1 \pmod{m}$.

We show that  $\sum_{i=1}^m x_i = m-1$ or $\ord g_1=m$. Assume the former does not hold, and let $\ell$ be maximal such that $\sum_{i=1}^{\ell}x_i= c<  m $.
We observe that $\s(k_0'^{m-\ell}(\prod_{i=1}^{\ell}f'_{-x_i})g_{1}^{c})= \s(k_0'^{m-\ell-1}(\prod_{i=1}^{\ell+1}f'_{-x_i})g_{1}^{[c+x_{\ell+1}]_m})$. By the choice of $\ell$ it follows that $[c+x_{\ell+1}]_m= c+x_{\ell+1}-m$.
Thus, $k_0'+ c g_1= k_{-x_{\ell +1}}+ (c+x_{\ell+1}-m)g_1$ and
$k_0'= k_{-x_{\ell +1}}+ (x_{\ell+1}-m)g_1$, which implies $mg_1=0$.

So, $\ord g_1= m$, or $S$ is of the form, where $k_0'=g_2$,
\[S= g_1^{sm-1}g_2^{(n-s)m}\prod_{i=1}^m ( - x_{i}g_1+g_2)\]
with $x_i \in [0, m-1]$ and $\sum_{i=1}^m x_i=m-1$.
Clearly, $\{g_1, g_2\}$ is a generating set of $G$.
Moreover, we know that if $s \ge 2$, then $\s(g_1^m)= e=\s(g_2^m)$.
Thus, $S$ is of the form given in (2).

Finally, suppose $\ord g_1= m$.
We have $s=1$. Let $\omega: G\to G / \langle g_1 \rangle$ denote the canonical map.
The sequence $\omega(\prod_{i=1}^{mn} h_i)$ is a mzss.  Thus, by Theorems \ref{thm_dir} and \ref{thm_invcyc}, $G / \langle g_1 \rangle$ is a cyclic group of order $mn$ and $\omega(\prod_{i=1}^{mn} h_i)= \omega(g_2)^{mn}$ for some $g_2\in G$, and $\ord g_2=mn$. Thus, $\{g_1,g_2\}$ is a basis of $G$ and $S$ has the form given in (1).

\medskip

Case 2: $\varphi(S)$ is of the form given in \eqref{eq_struc2}. If $m=2$, then $bf_1+2f_2= f_1$ and the sequence $\varphi(S)$ is also of the form given in \eqref{eq_struc1}. Thus, we assume $m \ge 3$. Moreover, we note that with respect to the basis  $f'_1 = f_1$ and $f'_2 = bf_1 + f_2$, we have $f_2= (m-b)f'_1+e'_2 $ and $bf_1 + 2f_2= (m-b)f'_1 + 2f'_2$. Thus, we may assume that $b<m/2$.
Let $S=FT$ with $\varphi(T)= f_1^{m}(bf_1 + f_2)^{ m-1}f_2^{m-1}(bf_1 + 2f_2)$.
We note that $F=\prod_{i=1}^{n-2}F_i$ with $\varphi(F_i) \in \{f_1^m, f_2^m, (bf_1 + f_2)^m\}$ for each $i \in [1, n-2]$.
Suppose $T=T_1T_2$ such that $\s(\varphi(T_i))=0$ and $T_i \neq 1$ for $i \in [1,2]$.
Then $\s(T_1)\s(T_2)\prod_{i=1}^{n-2}\s(F_i)$ is a ml-mzss over $H$ and thus equal to $e^n$ for some generating element $e$ of $H$.
It follows that for each  factorization $T=T_1T_2$ with $\s(\varphi(T_i))=0$ and $T_i \neq 1$ for $i \in [1,2]$, we have $\s(T_1)=\s(T_2)=e$.

Let $T_1'\mid T$ such that $\varphi(T_1')= f_1^b (bf_1 + f_2)^{ m-1} f_2$ and let $T_2'=T_1'^{-1}T$.
Suppose that, for some $i \in [1,2]$, there exists distinct elements $g_i,g_i'\in \supp(S)$ such that $\varphi(g_i)=\varphi(g_i')=f_i$.
We may assume that $g_i \mid T_1'$ and $g_i'\mid T_2'$.
It follows that $\s(g_i^{-1}g_i'T_1')=e=\s(T_1')$, a contradiction.
Thus, $\varphi^{-1}(f_i)\cap \supp(S)= \{g_i\}$ for $i \in [1,2]$.

Now, let $T_1'' \mid T$ such that $\varphi(T_1'')= f_1^{2b}(bf_1 + f_2)^{ m-2}f_2^2$ and $T_2''=T_1''^{-1}T$.
We can argue in the same way that  $\varphi^{-1}(bf_1+f_2)\cap \supp(S)= \{k_b\}$.
Finally, let $k\mid S$ such that $\varphi(k)=bf_1 + 2f_2$.
It follows that
\[ S = g_1^{s_1m}g_2^{s_2m-1} k_b^{s_3 m-1} k.\]
We note that $\ord g_1= mn$ and that $\{g_1, g_2\}$ is a generating set of $G$.

Since, as above, $\s(k_b^{m-2}kg_1^b)=e= \s(k_b^{m-1}kg_2^{m-1})$, it follows that  $bg_1=k_b +(m-1)g_2$.
Moreover, $\s(g_1^{2b}k_b^{m-2}g_2^2)= e= \s(g_1^{b}k_b^{m-1}g_2)$ implies that $bg_1+g_2= k_b$. Thus, $mg_2=0$ and and $\{g_2, g_1\}$ is a basis of $G$. Moreover, we have $s_2= 1$.
Additionally, we get $k+(m-1)g_2= bg_1+g_2$, i.e., $k =bg_1+2g_2 $.

We observe that the projection to $\langle g_1 \rangle$ (with respect to $G=\langle g_1 \rangle \oplus \langle g_2 \rangle$) of the sequence $g_2^{-(m-1)}S$ , i.e., the sequence $g_1^{s_1m}(bg_1)^{s_3m}$, is a mzss. Since $s_1+s_3=n$, this implies $b=1$.
Thus, $S$ is of the form given in (1).

\section*{Acknowledgment}
The author would like to thank the referee for suggestions and corrections, and A.~Geroldinger and D.~Grynkiewicz for valuable discussions related to this paper.


\end{document}